\documentstyle[amssymb,amsfonts]{amsart}

\def\beq{\begin{equation}}
\def\eeq{\end{equation}}
\def\barray{\begin{eqnarray*}}
\def\earray{\end{eqnarray*}}

\def\Vol{\hbox{\rm Vol}}
\def\var{\hbox{\bf Var}}
\def\Var{\hbox{\bf Var}}
\def\Conv{\hbox{ Conv}}
\def\R{{\hbox{\bf R}}}

\def\P{{\hbox{\bf P}}}
\def\pr{{\hbox{\bf P}}}
\def\E{{\hbox{\bf E}}}
\def\e{{\hbox{\bf E}}}

\font \roman = cmr10 at 10 true pt

\def\tri{{\triangle}}

\def\be#1{ \begin{equation}\label{#1} }

\def\bas{\begin{align*}}
\def\eas{\end{align*}}
\def\bi{\begin{itemize}}
\def\ei{\end{itemize}}

\def\aff{{\hbox{\roman aff}}}

\def \endprf{\hfill  {\vrule height6pt width6pt depth0pt}\medskip}
\def\emph#1{{\it #1}}
\def\textbf#1{{\bf #1}}

\def\BE{{\mathbf E}}

\def\BI{{\mathbf I}}

\def\CD{{\mathcal D}}

\def\BBR {{\mathbb R}}

\def\ep{{\epsilon}}

\def\hs{\hfill $\square$}

\theoremstyle{plain}
  \newtheorem{theorem}[subsection]{Theorem}

  \newtheorem{lemma}[subsection]{Lemma}
  \newtheorem{corollary}[subsection]{Corollary}
   \newtheorem{claim}[subsection]{Claim}

\theoremstyle{remark}
  \newtheorem{remark}[subsection]{Remark}

\theoremstyle{definition}

\include{psfig}

\begin{document}

\title[Gaussian polytopes]{Central limit theorems for Gaussian polytopes }

\author{ Imre B\'ar\'any}
\address{  R\'enyi Institute of Mathematics,  Hungarian Academy of Sciences,
 POBox 127, 1364 Budapest, Hungary  and Department of Mathematics, University College London,
    Gower Street, London WC1E 6BT, England.}
\email{barany@@renyi.hu}

\thanks{I. B\'ar\'any is supported by Hungarian National Foundation Grants T 046246
and T 037846}

\author{Van Vu}
\address{Department of Mathematics, Rutgers University, Piscataway, NJ 08854}
\email{vanvu@@math.rutgers.edu}

\thanks{V. Vu is an A. Sloan  Fellow and is supported by an NSF Career Grant.}

\begin{abstract} Choose  $n$ random, independent points in $\R^d$ according to the standard normal
distribution. Their convex hull $K_n$ is the {\sl Gaussian random
polytope}. We prove
that the volume and the number of faces of $K_n$ satisfy the central limit theorem, settling a
well known conjecture in the field.
\end{abstract}

\maketitle

\section{The main result}

Let $\Psi _d=\Psi $ denote the standard normal distribution on $\R^d$, its density
function is
\[
\psi _d=\psi = \frac 1 {(2\pi)^{d/2}}\exp\{-\frac {x^2}2\}
\]
where $x^2=|x|^2$ is the square of the Euclidean norm of $x \in \R^d$. We will use this
notation only for $d \geq 2$, for $d=1$ the standard normal has density function
\[
\phi = \frac 1 {(2\pi)^{1/2}}\exp\{-\frac {x^2}2\}
\]
with distribution $\Phi$.

Fix $d \ge 2$ and choose a set $X_n=\{x_1,\dots,x_n\}$ of random independent points from
$\R^d$ according to the normal distribution $\Psi$. The convex hull of these points,
$K_n=\Conv(x_1,\dots,x_n)$, is the {\sl Gaussian random polytope} or {\sl Gaussian
polytope} for short. This is one of the central models in the
 theory of random polytopes, initiated by R\'enyi and Sulanke in the 60s.
The main goal of this theory is to investigate the
 distributions of the key functionals (such as the volume) of random polytopes.

A cornerstone in probability theory is the central limit theorem.
 A sequence $X_n$ of random
variables satisfies the central limit theorem if for every $t$

$$ \lim_{n \rightarrow \infty} \P(\frac{X_n -\E X_n}{\sqrt {\Var X_n } } \le t) -\Phi (t) =0. $$

It is a natural and  important conjecture in the theory of random polytopes that the key
functionals of $K_n$ satisfy the central limit theorem, as $n$ tends to infinity. This
conjecture has been open  for several decades, and very few partial results have been
proved (see the next section).

In this paper, we develop a general frame work which enables us to confirm this
conjecture for many functionals.  Due to the length of the proofs,  we will focus on the
volume and the number of faces, perhaps the two most interesting parameters. Some other
functionals (such as the intrinsic volumes of the probability content) will be discussed
in Section 14.

For a convex polytope $K$, we use
$\Vol(K)$ and  $f_s(K)$ to denote its volume and number of faces of dimension $s$,
respectively. Here are our main results

\begin{theorem}\label{theo:CLT} Let $d$ be a fixed integer at least $2$.  There
is a function $\ep(n)$ tending to $0$ as $n$ tends to infinity such that the following
holds. For any value of $t$,
\begin{equation}  | \P \Big(\frac{\Vol (K_n) -
\E\Vol (K_n)}{ \sqrt {\Var \Vol (K_n)}} \le t \Big) -\Phi (t) | \le \ep(n).
\end{equation}
\end{theorem}

\begin{theorem}\label{theo:CLTfaces} Let $d$ be a fixed integer at least $2$ and $s$
be a non-negative integer at most $d-1$.  There is a function $\ep(n)$ tending to $0$ as
$n$ tends to infinity such that the following holds. For any value of $t$,

\begin{equation} \label{equ:CLT} | \P \Big(\frac{f_s (K_n) -
\E f_s (K_n) }{ \sqrt {\Var f_s (K_n)}} \le t \Big) -\Phi (t) | \le \ep(n).
\end{equation}
\end{theorem}

\begin{remark} \label{remark:rateofconv0} In both theorems,
we can take $\ep(n)= (\log n)^{-(d-1)/4+o(1)}$. (See Remarks \ref{remark:rateofconv1},
\ref{remark:rateofconv2} and \ref{remark:rateofconv3}.)
\end{remark}

In the next section, we  give a brief survey about the study
of Gaussian polytopes and random polytopes
in general.

{\bf Notation.} In the whole paper, we assume that $n$ is large, whenever needed. The
asymptotic notations are used under the assumption that $n \rightarrow \infty$. Given
non-negative functions $f(n)$ and $g(n)$, we write $f(n)=O(g(n))$ ($f(n) =\Omega (g(n))$)
if there is a positive constant $C$, independent of $n$, such that $f(n) \le Cg(n)$
($f(n) \ge C g(n)$)  for all sufficiently large value of $n$. We write $f(n) =\Theta
(g(n))$ if $f(n)= O(g(n))$ and $f(n)=\Omega (g(n))$. In this case, we say that $f(n)$ and
$g(n)$ have the same order of magnitude. Finally $f(n) =o(g(n))$ if $f(n)/g(n)$ tends to
zero as $n$ tends to infinity.

Consider a  (measurable) subset $S$ of $\R^d$. The probability content of $S$ is
$$\Psi (S) = \int_S \psi(x) d x. $$
$\P$, $\E$, $\Var$ denote probability, expectation, variance, respectively. Let $t_i $,
$i=1, \dots, n$, be independent random variables and $Y=Y(t_1, \dots, t_n)$ be a random
variable depending on $t_1, \dots, t_n$. $\BE(Y|t_1, \dots, t_i)$ is the conditional
expectation of $Y$ conditioned on the first $i$ variables. ${\BI}_E$ is the indicator of
the event $E$: $\BI_E=1$ if $E$ holds and $0$ otherwise.

\section{History}

Gaussian random polytopes were first considered by R\'enyi and Sulanke in their classical
paper \cite{rs1}. Naturally, the existence of central limit theorems should be one of the
very first questions to ask. However, early results are very far from a possible answer
of this question, due to the lack of tools. These results
  mostly focused on expectations. In particular, R\'enyi and Sulanke
determined the expectation of $f_1(K_n)$ for a Gaussian polytope in $\R^2$. (Here and
later $f_i$ denotes the number of faces of dimension $i$.) In 1970,  Raynaud \cite{ray}
computed $\e f_{d-1}(K_n)$ in all dimensions. The general formula is
\begin{equation}
\e f_s(K_n)=\frac {2^d}{\sqrt d}{d \choose {s+1}}\beta_{s,d-1}(\pi \log n)^{\frac
{d-1}{2}}(1+o(1)) \label{eq:fs}
\end{equation}
where $s \in \{0,1,\dots,d-1\}$ and $d \geq 1$, as $n \rightarrow  \infty$. Here $\beta_{s,d-1}$
is the internal angle of the regular $(d-1)$-simplex at one of its $s$-dimensional faces.
The formula was proved by Affentranger and Schneider~\cite{as} and by Baryshnikov and
Vitale~\cite{bv}; simpler proofs can be found in \cite{hmr}.  Recently Hug
and Reitzner~\cite{hr} obtained an estimate for the variance

 \begin{equation}
\var f_s(K_n) =O( (\log n)^{\frac {d-1}2}). \label{eq:varfs}
\end{equation}

In \cite{hue1, hue2}, Hueter stated a central limit theorem for $f_0(K_n)$, but the proof had a
gap, namely, the claimed estimate on the variance was not correct.

As far as the volume is concerned, Affentranger~\cite{aff} determined the expectation of
$\Vol (K_n)$:
\begin{equation}
\e \Vol (K_n)=\kappa_d (2\log n)^{\frac d2}(1+o(1)). \label{eq:vol}
\end{equation}
Here $\kappa_d$ denotes the volume of $B^d$, the $d$-dimensional unit ball. An upper
bound for the variance of $\Vol (K_n)$ is given by Hug and Reitzner~\cite{hr}:
\begin{equation}
\var \Vol(K_n) =O( (\log n)^{\frac {d-3}2})  \label{eq:varvol}
\end{equation}

We are not aware of a central limit theorem for the volume, prior to this paper.

Another popular model of random polytopes is the so-called uniform model, defined as
follows. Let $K$ be a convex set in $\R^d$ of volume one. Select $n$ random points in $K$
with respect to the uniform distribution and define the random polytope as the convex
hull of these points. Similar to the situation with the Gaussian model, there is a vast
amount of literature focusing on the expectations of the key functions (see \cite{ww} for
a survey). As far as central limit theorems are concerned, the case $d=2$ has been
studied  by Groeneboom~\cite{groe}, Groeneboom and Cabo~\cite{cg}, and Hsing~\cite{hs}.
They proved central limit theorems for random polyogon in the square and the unit disk.
But their methods do not extend to higher dimensions.

In 2004 and 2005 there were several notable developments on the uniform model, especially
in the case when the mother body $K$ has  smooth boundary:  Vu \cite{vvu1} proved that
several key functionals have distributions with exponential tails. Next,
Reitzner~\cite{rei1} established a central limit theorem for a Poisson variant of the
model. Further, Vu \cite{vvu2}, using the results of the above two papers and a coupling
argument, proved several central limit theorems for the uniform model. The central limit
theorem when $K$ is a polytope was established by B\'ar\'any and Reitzner \cite{br}.

\vskip2mm

The frame work we develop in this paper makes use of  ideas  from \cite{rei1, vvu1, vvu2}
and also from \cite{br}. Moreover, due to the obvious differences between the  uniform measure
and  the Gaussian one, we also need to introduce several new ideas to handle technical obstacles.

\vskip2mm

Let us conclude this section with a few basic facts about the normal distribution. Let
$r$ be a positive number at least one. Let $B(r)$ denote the ball of radius $r$ centered
at the origin and $\overline {B(r)}$ be its complement. The probability content of
$\overline {B(r)}$ is

\begin{equation} \label{equ:contentofball} \Psi( \overline {B(r)} )  =\Theta ( e^{-r^2/2} r^{d-2}).
\end{equation}

Let $H(r)$ be a half space at distance $r$ from the origin ($H(r)$ is not unique, but it
does not matter). The probability content of $H(r)$ is

\begin{equation} \label{equ:contentofspace} \Psi( H(r)  )  =\Theta ( e^{-r^2/2} r^{-1}).
\end{equation}

\section{ Two more models} \label{section:2models}

It is hard to prove the CLT for $K_n$ directly. We are going to take a detour and prove the
CLT for some more convenient models, namely $K'_n$ and $\Pi_n$, and next prove that the
distributions of $\Vol (K_n)$ and $\Vol (K'_n)$ and $\Vol (\Pi_n)$ are approximately the
same.

We define $K'_n$ first. Let $c_0$ be a large constant compared to the dimension $d$
($c_0=100d$ will satisfy all purposes).  Define $R>0$ via
\begin{equation} \label{equ:defR} R^2=2\log n+ \log(\log n)^{c_0}.
\end{equation}
We will use this definition later as well, for the time being we only need the following
consequence.
\begin{equation} \label{equ:defR0} e^{-R ^2/2} R ^{d-2}  = \Theta \left(\frac{(\log n)^{(d-2)/2}}{n (\log n)^{c_0/2}}\right)
= \Theta \left(\frac{1}{n (\log n)^{C_0}}\right).
\end{equation}
where $C_0=\frac {c_0}2-\frac {d-2}2$. Notice that the left hand side is (up to a
constant factor) the probability content of the complement of $B(R)$,  the ball of radius
$R$ centered at the origin, see \eqref{equ:contentofball}. The probability that one of
$n$ random points falls outside $B(R)$ is at most

$$ O( n \times \frac{1}{n (\log n)^{C_0}}) = O( \frac{1}{ (\log n)^{C_0}}). $$

By setting $c_0$ (and so $C_0$) sufficiently large, this probability will be negligible.
This allows us to replace the normal distribution $\Psi$ by the truncated distribution
$\Psi'$, restricted to $B(R)$. $\Psi'$ is defined so that for any region $S$ in $B(R)$,
the measure of $S$ is $\Psi'(S)= \frac{ \Psi(S)}{\Psi(B(R)}$. To be precise, the density
function $\psi'$ of $\Psi'$ is defined as

$$ \psi' (x) =\psi (x)  \frac{\BI_{x \in B(R)} }{\Psi(B(R))} $$

\noindent where $\BI$ is the indicator variable.

Let $K_n'$ be the convex hull of a set of $n$ random points chosen independently in
$B(R)$ with respect to $\Psi'$. The central limit theorem for the $K'_n$ model says the
following.

\begin{theorem} \label{theo:CLT1} Let $d$ be a fixed integer at least $2$. There is a
function $\ep(n)$ tending to zero as $n$ tends to infinity such that

$$ | \P \Big( \frac{\Vol(K'_n))- \E(\Vol(K'_n))} {\sqrt{\Var \Vol(K'_n)}}  \le t \Big)-\Phi (t) | \le \ep(n) $$

\noindent holds for all $t$.

\end{theorem}

Again, it is hard to prove this theorem directly. That's why we need the second model,
the Poisson polytope.

We consider a Poisson point process, $X(n)$, of intensity $n$ and underlying distribution
$\Psi'$ where $\Psi'$ is the truncated Gaussian, that is, the Gaussian restricted to
$B(R)$. Let $S$ be a measurable subset of $\R^d$. The intersection of $X(n)$ with $S$
consists of random points $\{x_1,\dots,x_k\}=X(n) \cap S$ where the number, $k$, of
random points is Poisson distributed with expectation $n\Psi'(S)$ and for fixed $k$, the
points are distributed independently. The property that we need most is that if $S_1$ and
$S_2$ are disjoint measurable sets, then the two point sets
$\{x_1,\dots,x_{k_1}\}=X(n)\cap S_1$ and $\{y_1,\dots,y_{k_2}\}=X(n)\cap S_2$ are
independent, $k_1$ and $k_2$ are independently Poisson distributed. The {\sl Poisson
polytope} is, by definition, the convex hull of $X(n)$.

Another, equivalent and useful, way to look at $\Pi_n$ is the following. First choose a
random number $n'$ with respect to the Poisson distribution with mean $n$. Next, generate
$n'$ random, independent points $x_1,\dots,x_{n'}$ with respect to $\Psi '$, the
truncated normal distribution on $\R^d$. Then  $\Pi_n$ is the convex hull,
$\Conv\{x_1,\dots,x_{n'}\}$, of the chosen points. It is well known that $n'$ is very
close to $n$ with high probability:
\[
\P(|n'-n| \geq A\sqrt {n\log n})\leq n^{-A/4},
\]
for every constant $A\ge 10$ (the constants $4$ and $10$ are just convenient choices and
play no important role). So a good approximation of the Poisson polytope $\Pi_n$ is
$K_{n'}$ with $n'$ Poisson distributed. Clearly, $n'$ is concentrated on the interval
$I=[n-A\sqrt {n\log n},n+A\sqrt {n\log n}]$ and negligible outside this interval. The
central limit theorem for the Poisson model is as follows.

\begin{theorem} \label{theo:Poisson1}
Let $d$ be a fixed integer at least $2$.  There is a function $\ep(n)$ tending to $0$ as
$n$ tends to infinity such that the following holds. For any value of $t$,

$$ | \P \Big(\frac{|\Vol (\Pi_n) - \E\Vol (\Pi_n)| }{ \sqrt {\Var \Vol (\Pi_n)}} \le t
\Big) -\Phi (t) | \le \ep(n).
$$

\end{theorem}

\begin{remark} \label{remark:rateofconv2} In both theorems above one can take
$\ep(n)= (\log n)^{-(d-1)/4+o(1)}$. This error term will be the dominating one when we
apply Lemma \ref{lemma:twodistribution} from the next section.
\end{remark}

\section{ The plan  of the proof} \label{section:plan}

From now on we focus on the volume, the proof for the number of faces is basically
the same and will be discussed in Section \ref{section:faces}.

The proof is long and consists of many steps. To help the reader grasp the main ideas
quickly, we first lay out the plan of the proof. The leading idea is coupling. In fact,
our proof will involve two different couplings. Both of them are based on a simple lemma.

\begin{lemma} \label{lemma:twodistribution} Let $Y_n$ and $Y_n'$ be two sequences of random variables
with means $\mu_n$ and $\mu'_n$, variances $\sigma_n^2$ and ${\sigma'}_n^2$,
respectively. Assume that there are functions $\ep_1(n), \ep_2(n)$, $\ep_3(n), \ep_4(n)
$, all tending to zero as $n$ tends to infinity such that

\begin{itemize}

\item $ |\mu'_n -\mu_n| \le \ep_1 (n) \sigma'_n $

\item $ | {\sigma'}_n  - \sigma_n| \le \ep_2(n) {\sigma'}_n $

\item For any $t$, $|\P(Y'_n \ge t) - \P(Y_n \ge t)|  \le \ep_3 (n) $

\item  For any $t$,
$$| \P \Big(\frac{Y'_n- \mu'_n}{ \sigma'_n} \le t\Big) - \Phi (t) | \le \ep_4(n) .$$
\end{itemize}

Then there is a positive constant $C$  such that
 for any $t$,

$$| \P \Big(\frac{Y_n- \mu_n}{ \sigma_n} \le t\Big) - \Phi (t) | \le C\sum_{i=1}^4 \ep_i(n) .$$

\end{lemma}

Basically, this lemma asserts that if $Y'_n$ satisfies the CLT (the fourth condition) and
$Y_n$ is sufficiently close to $Y'_n$ in distribution (the first three conditions), then
$Y_n$ also satisfies the CLT. We defer the routine proof to the end of this section. The
lemma has been used in an implicit form in \cite{vvu2} and in \cite{rei1}.

\begin{remark} \label{remark:rateofconv1} We can rewrite the error term
$C\sum_{i=1}^4 \ep_i(n) $ as $C \max_{i=1}^4 \ep_i (n)$ (the two $C$'s can have different
values). In applications of Lemma \ref{lemma:twodistribution}, $\ep_4(n)$ will be the
dominating term.
\end{remark}

\vskip2mm We now present the plan for the proof of  Theorem \ref{theo:CLT}, which
consists of the following  steps.
\begin{itemize}

\item {\bf Step 1.} (Variance) In this step,  we show that the exact order of magnitude
of $\Var \Vol(K_n)$ is $(\log n)^{(d-3)/2}$. The upper bound was obtained in \cite{hr}.
We will prove the matching lower bound. Section \ref{section:step1} is devoted to this
step. The necessary geometric tools are developed in Section \ref{section:geom}. The
variance plays a significant role and we will use the estimate obtained in this step
several times later on.

\vskip2mm

\item {\bf Step 2.} (The first coupling) In this step, we couple $K_n$ and $K'_n$ in
order to show that they satisfy the first three conditions of Lemma
\ref{lemma:twodistribution}. This will be done in Section \ref{section:couple1}. Thus, it
remains to verify the fourth, and critical, condition that $\Vol (K'_n)$ satisfies the
CLT. This task will take time and effort. We mention that the second condition of Lemma
\ref{lemma:twodistribution}, together with Step 1, imply that the order of magnitude of
$\Var \Vol(K'_n)$ is $(\log n)^{(d-3)/2}$.

\vskip2mm

\item {\bf Step 3.} (The second coupling) In this step which is in Section \ref{section:couple2},
we couple $\Pi_n$ with $K'_n$. Technically speaking, we are going to verify the first
three conditions of Lemma \ref{lemma:twodistribution} with respect to $\Vol (\Pi_n)$ and
$\Vol (K'_n)$. After this, both Theorem \ref{theo:CLT} and Theorem \ref{theo:CLT1} follow
from Theorem \ref{theo:Poisson1}, the CLT for the Poisson model. This step is close to
the coupling argument used for the uniform model \cite{vvu2}. However, the analysis for
the current case is simpler, as strong concentration results are not needed. Again, the
results imply that the order of magnitude of $\Var \Vol(\Pi_n)$ is $(\log n)^{(d-3)/2)}$.

\vskip2mm

\item {\bf Step 4.} (Sandwiching)  In this step, we define a radius $r < R$ but very
close to $R$, and prove that $K'_n$ contains the ball $B(r)$ with high probability,
namely, with probability $1-(\log n)^{-C}$. (For this end $r$ has to be chosen carefully,
see Remark 9.4.) By definition, $K'_n$ is contained in $B(R)$. So with high probability,
$K'_n$ is sandwiched between two very close balls. We will also prove that the Poisson
polytope has the same property, that is, $B(r) \subset \Pi_n \subset B(R)$ with high
probability. This is the content of Section \ref{section:sandwich}.

\end{itemize}

The main idea behind the proof of Theorem \ref{theo:Poisson1}, following Reitzner
\cite{rei1}, is as follows. It is well known that if $\xi_1, \dots, \xi_n$ are
independent variables with bounded means and variances, then the distribution of the
normalized version of the sum $\sum_{i=1}^n \xi_i$ is approximately Gaussian. We are
going to use a strengthening of this result, originally due to C. Stein \cite{st}, which
asserts that it suffices to assume that the $\xi_i$ are weakly dependent. The
quantitative, and technical, statement below is from Rinott \cite{ri}, which is  slightly
stronger than an earlier one due to Baldi and Rinott \cite{BaRi}.

\begin{theorem}\label{th:rinott} Assume $G$ is a graph with vertex set $V(G)$ and edge set $E(G)$,
$|V(G)|=m$, and maximal degree $D$. Assume $\xi_v$ is a random variable satisfying
$|\xi_v|\leq M$ almost surely for each $v \in V(G)$. Assume further that if there is no
edge between a vertex in $V_1 \subset V(G)$ and a vertex of $V_2 \subset V(G)$ where
$V_1$ and $V_2$ are disjoint, then the random variables $\{\xi_v:\; v \in V_1\}$ and
$\{\xi_v:\; v \in V_2\}$ are independent. Then, writing $\xi=\sum_{v \in V(G)} \xi_v$, we
have
\[
\left|\pr\left(\frac {\xi-\e \xi}{\sqrt{\var \xi}} -\Phi(t)\right)\right| \leq \frac
{DM}{\sqrt {\var \xi}}\left(\frac {1}{\sqrt{2\pi}}+ 16\frac {\sqrt {mD}M}{\sqrt {\var
\xi}}+10 \frac {mDM^2}{\var \xi}\right).
\]
\end{theorem}

In order to apply this result we have to make some geometric preparations and define the
{\sl dependency graph}.

\begin{itemize}

\item {\bf Step 5.} (The dependency graph) We subdivide the annulus $A(R,r)= B(R)
\setminus B(r)$ into pairwise internally disjoint cells $W_1,\dots,W_m$. The cells are
nice and well-behaving, and they define the {\sl dependency graph} $G$ with vertex set
$V(G)=\{1,\dots,m\}$ and the pair $(i,j)$ forming an edge of $G$ if $W_i$ and $W_j$ are
far apart. (The actual definition is different, but this is the essence of it.) Note that
the dependency graph is defined by geometric conditions. We will give an upper bound on
the maximal degree of $G$, and on the volume of the cells. The details appear in Section
~\ref{section:depend}. Note that randomness does not come up here but is present in the
background.

\vskip2mm

\item {\bf Step 6.} (CLT for the Poisson model) In this step, we work with the
Poisson model $\Pi_n$ under condition $B$ which says that $B(r) \subset \Pi_n$. The
Baldi-Rinott theorem can be applied with $\xi_i=\Vol (\Pi_n \cap W_i)$ and dependency
graph $G$. This is a technical step which is carried out in Section
\ref{section:CLTPoisson}. It proves Theorem \ref{theo:Poisson1}, the CLT for the Poisson
model, but only under condition $B$. The role of the Poisson model is critical here, as
it guarantees that $\xi_i$ and $\xi_j$ are independent whenever $i$ and $j$ are not
adjacent in $G$.

\vskip2mm

\item {\bf Step 7.} (Removing condition $B$) This is a technical step which is another,
(this time simple) application of Lemma~\ref{lemma:twodistribution}. It proves, finally,
that $\Vol (\Pi_n)$ satisfies the CLT (Theorem \ref{theo:Poisson1}) and so it finishes
the proof of the main theorem.

\end{itemize}

The proof for Theorem \ref{theo:CLTfaces} concerning the number of faces is similar and
will be presented in Section \ref{section:faces}. In the last  Section \ref{section:remarks}, we discuss few
 other results which can be proved using the same
method.

Let us now conclude this section with the proof of Lemma \ref{lemma:twodistribution}.

{\bf \noindent Proof of Lemma \ref{lemma:twodistribution}.} We have to show that for any
$x$

$$\P \Big( \frac{Y_n -\mu_n}{\sigma_n} \le x\Big) = \Phi (x) + O(\sum_{i=1}^4 \ep_i (n)). $$

\noindent By the third condition of the lemma

$$ \P ( \frac{Y_n -\mu_n}{\sigma_n} \le x) = \P( Y_n \le \mu_n + x \sigma_n)
= \P( Y'_n  \le \mu_n + x \sigma_n) + O(\ep_3(n)) . $$

\noindent On the other hand,

$$ \P( Y'_n \le \mu_n + x \sigma_n) = \P (Y'_n \le \mu'_n + x' \sigma'_n) $$

\noindent where $x'= \frac{\mu_n -\mu'_n}{\sigma'_n} + \frac{x \sigma_n}{\sigma'_n} $.
 The first two conditions of the lemma guarantee that
 $x'$ is between the maximum and minimum of the four values
 $   x (1 \pm \ep_2(n)) \pm \ep_1(n).$ Moreover,
the fourth condition of the lemma yields

$$ \P (Y'_n \le \mu'_n + x' \sigma'_n ) = \Phi (x') + O(\ep_4(n)). $$

\noindent Further,

$$\Phi (x') = \Phi (x) + (x'-x) \Phi' (x_0) $$

\noindent for some $x_0$ between $x$ and $x'$. The difference $|x-x'|$ is
at most  $|x| \ep_2(n) + \ep_1(n) $. As $\Phi' (x)$ decays
exponentially, it is easy to see that
$|x| \Phi'(x_0) = O(1)$ and thus

$$\Phi(x') = \Phi (x) + O(\ep_1(n) + \ep_2(n)). $$

\noindent Putting everything together completes the proof:

$$\P\Big(\frac{Y_n -\mu_n}{\sigma_n} \le x\Big) = \Phi (x) + O(\ep_1(n)+\ep_2(n)+\ep_3(n) +\ep_4(n))
. $$ \hs

\vskip3mm

\section{A geometric construction} \label{section:geom}

Here we give a geometric construction, \'a la Reitzner~\cite{rei1} and B\'ar\'any,
Reitzner~\cite{br}. We use it in the next section for estimating $\var \Vol(K_n)$ and
$\var f_s(K_n)$. A similar, if more subtle, construction will be needed for the
dependency graph as well.

In the construction $b_1,b_2,\dots$ are positive constants that depend on dimension only.
Let $S(r)$ denote the sphere of radius $r$ centered at the origin. We define
\[
r^2=2\log n- \log \log n.
\]
The choice of $r$ is not arbitrary here: it ensures that $\Psi(\tri_i)=\Theta(1/n)$ (see
later). Next we choose a system of points $y_1,\dots,y_m$ from the sphere $S(r)$ which is
maximal with respect to the property that for distinct $i$ and $j$
\[
|y_i-y_j| \geq 2b_1.
\]
Such a system can be found by an obvious greedy algorithm. The spherical caps on $S(r)$
with centre at $y_i$ and radius $b_1$ are pairwise disjoint, and the same spherical caps
with radius $2b_1$ cover $S(r)$. This implies by volume comparison
\begin{claim}\label{cl:size}
\[
m = \Theta\left((\log n)^{\frac {d-1}2}\right).
\]
\end{claim}

Next, for each $i=1,\dots,m$ set
\[
y_i^0=\left(1+\frac {1}{r^2}\right)y_i.
\]
Thus $|y_i^0|=r+\frac 1r$ and we have, for all $x \in \R^d$ with $r\leq |x| \leq r+\frac
1r$ that
\begin{equation}
\psi(x)=\Theta \left(\frac {\sqrt {\log n}}{n}\right) \label{eq:psi}
\end{equation}

Next we let $H_i$ denote the hyperplane with equation $z\cdot y_i =r^2$. For each
$i=1,\dots,m$ we fix a regular $(d-1)$-dimensional simplex in $H_i$ whose vertices
$y_i^1,\dots,y_i^d$ lie in the $(d-2)$-dimensional sphere
\[
H_i \cap S(y_i,\sqrt 2).
\]
The centre of this simplex is clearly $y_i$. The simplex $\tri_i$ is now defined as the
convex hull of the $y_i^j$, $j=0,1,\dots,d$.

\begin{claim}\label{cl:tri} For all $i$
\[
\Psi(\tri_i) =\Theta\left( \frac 1n \right).
\]
\end{claim}

{\bf \noindent Proof.} It is clear that for $j=1,\dots,d$
\[
|y_i^j| = \sqrt {r^2+2}< r+\frac 1r =|y_i^0|.
\]
Then every $x \in \tri _i$ satisfies $r\leq |x| \leq r+\frac 1r$, and the claim follows
from (\ref{eq:psi}) as $ \Vol \tri_i = \Theta\left(\frac 1{\sqrt {\log n}}\right)$. \hs

As the final step of the construction, for $i=1,\dots,m,\; j=0,1,\dots,d$,  let
$\tri_i^j$ be a homothetic copy of $\tri _i$ where the centre of homothety is $y_i^j$ and
the factor of homothety is a small number $b_2>0$.

This is our geometric construction. Now we establish several properties of this
construction.

\begin{claim}\label{cl:smalltri}
\[
\Psi(\tri_i ^j) =\Theta\left( \frac 1n \right).
\]
\end{claim}

{\bf \noindent Proof.} The density $\psi(x)$ satisfies (\ref{eq:psi}) for all $x \in \tri
_i^j$. The claim follows as the volume of $\tri_i^j$ is just $b_2^d$ times that of $\tri
_i$. \hs

Assume now that $z_j$ is an arbitrary point in $\tri_i^j$, $j=0,1,\dots,d$. We define the
cone $C_i$ via
\[
C_i=z_0+\mbox{pos}\{z_j-z_0:\; j=1,\dots,d\}.
\]

The following lemma is crucial since it implies the independence structure of $K_n$
needed when estimating the variance.

\begin{lemma}\label{le:cone} For $b_1$ large enough and $b_2$ small enough the cone
$C_i$ contains all simplices $\tri _k$ with $k\ne i$.
\end{lemma}

{\bf \noindent Proof.} We have to check that the segment $[z_0,y_j^k]$ intersects $\Conv
\{z_1,\dots,z_d\}$ whenever $j\ne i$ and $k\in \{0,1,\dots,d\}$. This is the same as
checking that the segment $[z_0,y_j^k]$ intersects $\Conv \{z_1',\dots,z_d'\}$ where
$z_j'=\aff\{z_0,z_j\}\cap H_i$. If $b_2$ is small enough then the $(d-1)$-dimensional
ball $B_i=H_i \cap B(y_i,\frac {\sqrt 2}{2d})$ is contained in $\Conv
\{z_1',\dots,z_d'\}$. It is not hard to see that, for large enough $b_1$, the segment
$[y_i^0,y_j^k]$ intersects $H_i \cap B(y_i,\frac {\sqrt 2}{3d})$ which is a smaller
shrunken copy of $B_i$.  (Here again $j\ne i$ and $k\in \{0,1,\dots,d\}$.) But $z_0$ is
very close to $y_i^0$ if the factor of homothety, $b_2$ is very small, and then the
segment $[z_0,y_j^k]$ intersects $B_i$. \hs

We need one more lemma for estimating the variance. Let $H_i^j$ be the halfspace
containing $\tri_i^k$ for all $k=1,\dots,d$ except $k= j$, not containing $\tri_i^0$ and
$\tri_i^j$, and whose bounding hyperplane touches all $\tri_i^k$ except $k=j$.

\begin{claim}\label{cl:measH} If $b_2$ is small enough, then
\[
\Psi(H_i^j) =O(n^{-1}).
\]
\end{claim}

{\bf \noindent Proof.} Let $H$ denote the hyperplane through the points $y_i^k$
($k=0,1,\dots,d$, $k\ne j$) for this proof. It is not hard to check that the distance of
$H$ from the origin is at least $r-\frac {d^2}r$. The bounding hyperplane of $H_i^j$
tends to $H$ as $b_2$ tends to zero. So for small enough $b_2$, the distance of $H_i^j$
from the origin is at least $r-\frac {2d^2}r$. An application of
(\ref{equ:contentofspace}) finishes the proof.\hs

\section{The variance} \label{section:step1}
\begin{theorem}\label{th:varvol} $\Var \Vol (K_n) =\Theta ( (\log n)^{\frac {d-3}2})$.
\end{theorem}

{\bf \noindent Proof.} The upper bound (\ref{eq:varvol}) has been proved by Hug and
Reitzner~\cite{hr}. So we need to give a lower bound on $\Var \Vol (K_n)$.

Let $X_n=\{x_1,\dots,x_n\}$ denote our random sample of $n$ points. Denote by $A_i$ the
event that exactly one random point (out of the sample $X_n$) is contained in each
simplex $\tri_i^j$, $j=0,1,\dots,d$ and no further point of $X_n$ is contained in $H_i^+
\cup \bigcup_{j=1}^d H_i^j$. Here $H_i^+$ is the halfspace not containing the origin
whose bounding hyperplane is $H_i$. Since $H_i^+$ is farther from the origin than $H_i^j$
($j>0$), Claim~\ref{cl:measH} implies $\Psi(H_i^+) =O(1/n)$.

\begin{lemma}\label{le:condA} There is a positive constant $b_3$ such that, for every $i=1,\dots,m$
\[
\P(A_i) \geq b_3.
\]
\end{lemma}

{\bf \noindent Proof.} Assuming that $A_i$ has occurred, let $x_j \in X_n$ denote the
unique point of $X_n$ in $\tri_i^j$, $j=0,1,\dots,d$, and set $X=X_n\setminus
\{x_0,\dots,x_d\}$. As $\Psi (\tri _i^j)=\Omega (1/n)$ and $\Psi (H_i^j) =O( 1/n)$ we
have
\begin{eqnarray*}
\P(A_i)&=& {n \choose {d+1}} \P(x_j \in \tri_i^j,\; j=0,\dots,d)\P(X \cap (H_i^+ \cup \bigcup_{k=1}^d H_i^k)=\emptyset)\\
          &=& {n \choose {d+1}} \prod_0^d \Psi (\tri_i^j)\left(1-\Psi(H_i^+              \cup \bigcup_{k=1}^d H_i^k)\right)^{n-d-1}\\
          &\ge& c_1 n^{d+1}\cdot \frac 1{n^{d+1}}\left(1-\frac cn\right)^{n-d-1}\ge b_3>0.
\end{eqnarray*}
Here $c$ is $(d+1)$ times the implicit constant in Claim~\ref{cl:measH}, and $c_1$ is
another constant that depends on $d$ only. \hs

So we can bound the expected number of $A_i$ from below:
\[
\e\left(\sum_1^m \BI _{A_i}\right)=\sum_1^m \P(A_i)=\Omega (m). \label{eq:expA}
\]

We start bounding $\var \Vol (K_n)$ from below. Let $\cal F$ denote the position of all
random points from $X_n$ except those in $\tri_i^0$ with $\BI_{A_i}=1$, $i=1,\dots,m$. We
decompose the variance under condition $\cal F$:
\begin{equation}
\var \Vol(K_n)= \e \var (\Vol(K_n)|{\cal F})+\var \e (\Vol(K_n)|{\cal F}) \geq \e \var
(\Vol(K_n)|{\cal F}). \label{eq:decomp}
\end{equation}
Suppose condition $\cal F$ holds and ${\BI}_{A_i}={\BI}_{A_j}=1$. Clearly,  the unique
$x_i \in \tri_i^0$ and $x_j \in \tri_j^0$ ($x_i,x_j \in X_n$) are vertices of $K_n$, and,
because of Lemma~\ref{le:cone}, there is no edge between $x_i$ and $x_j$. Then the change
in $K_n$ when $x_i$ is moved is independent of the change when $x_j$ is moved. This
implies that the change in $\Vol(K_n)$ when $x_i$ is moved is independent of the change
when $x_j$ is moved, showing that
\[
\var (\Vol(K_n)|{\cal F}) = \sum_{i:\; \BI_{A_i}=1} \var _{x_i}\Vol(K_n)
\label{eq:sumvar}
\]
where the variance in the sum is taken when $x_i$ is changing within $\tri _i^0$.

We now evaluate this variance. Let $z_j \in X_n$ be the unique random point in
$\tri_i^j$, ($j=1,\dots,d$). Denote the simplex $\Conv \{x_i,z_1,\dots,z_d\}$ by
$\triangle$. The change in $\Vol(K_n)$ when $x_i$ changes within $\tri_i^0$ equals the
change in $\Vol(\triangle)$ and
\[
\var _{x_i}\Vol(\tri) = \e\left(\Vol(\tri) - \e_{x_i}\Vol(\tri)\right)^2.
\]
The base of $\tri$, $\Conv\{z_1,\dots,z_d\}$, is a fixed $(d-1)$-dimensional simplex, of
constant $(d-1)$-dimensional volume. Its height varies nearly between $\frac 1r (1-b_2)$
and $\frac 1r$, so the expectation $\e_{x_i}\Vol(\tri)$ is about $\Theta(1/r)$. Moreover,
the height of $\tri$ changes on a small interval of length about $b_2/r$, so the volume
is a linear (but not constant) function on a positive fraction of this interval.
Consequently,
\[
\left(\Vol(\tri) - \e_{x_i}\Vol(\tri)\right)^2 =\Omega\left(\frac 1{(\sqrt {\log
n})^2}\right) = \Omega\left( \frac 1{\log n}\right)
\]
holds on a positive fraction of $\tri _i^0$. This implies that
\[
\var _{x_i}\Vol(\tri) =\Omega \left( \frac 1{\log n}\right).
\]
Putting this into formula (\ref{eq:decomp}) and using (\ref{eq:expA}) completes the
proof. \hs

The same method, with the same notation, works for $\var f_s(K_n)$, so we present it
here.

\begin{theorem}\label{th:varfs} $\var f_s (K_n) =\Theta\left( (\log n)^{\frac {d-3}2}\right)$
\end{theorem}
{\bf \noindent Proof.} The upper bound is again due to Hug and Reitzner~\cite{hr}.

The method for the lower bound is similar to the one in Reitzner~\cite{rei1}. We assume
$s\in \{0,1,\dots,d-1\}$. Condition $A_i$ is the same as in Lemma~\ref{le:condA} except
that we require exactly two points from $X_n$ to be in $\tri _i^0$. Also, we let $\cal F$
denote the position of all random points from $X_n$ except those two in $\tri _i^0$ with
$\BI_{A_i}=1$, $i=1,\dots,m$. Then Lemma~\ref{le:condA} remains valid for the new $A_i$.
We can decompose the variance under condition $\cal F$ the same way and we still get
(\ref{eq:decomp}). An identical analysis applies and gives
\[
\var (f_s(K_n)|{\cal F}) \ge \sum_{i\: \BI_{A_i}=1} \var _{x_i,y_i}f_s(K_n)
\]
where the variance in the sum is taken when $x_i,y_i$ are changing within $\tri _i^0$.
Here $x_i$ and $y_i$ are the two points from $X_n$ contained in $\tri_i^0$. The proof of
the following claim is simple and left as an exercise.
\begin{claim}
\[ \var _{x_i,y_i}f_s(K_n) =\Theta(1).
\]
\end{claim}
This finishes the proof of Theorem~\ref{th:varfs}. \hs

\section{The first coupling} \label{section:couple1}

Here we show that the random variables $\Vol (K_n)$ and $\Vol (K'_n)$ satisfy the first
three conditions of Lemma \ref{lemma:twodistribution}.

\begin{lemma} \label{lemma:coupling1} We have

$$| \E\Vol (K'_n)- \E\Vol (K_n)| \le \sqrt {\Var \Vol(K_n)} (\log n)^{-C_0/2} $$

$$| \Var \Vol(K'_n)- \Var \Vol(K_n)| \le  \Var \Vol(K_n) (\log n)^{-C_0/2}. $$

\noindent Furthermore, for all $t$,

$$|\P (\Vol (K'_n) \ge t) - \P (\Vol (K_n) \ge t) | \le (\log n)^{-C_0/2}. $$
\end{lemma}

{\bf \noindent Proof of Lemma \ref{lemma:coupling1}.} Choose $n$ points $t_1, \dots, t_n$
in $\R^d$ with respect to the normal distribution $\Psi$. Let $A$ denote the event that
all $n$ points fall inside $B(R)$. (Recall that $R$ is defined in (\ref {equ:defR}).) For
every non-negative integer $i$, let $B_i$ be the event that all $n$ points fall inside
$B(4^{i+1} R)$ but there is at least one point outside $B(4^ i R)$. Trivially

$$ \overline {A} = \cup_{i=0}^{\infty} B_i. $$

Let $Y=Y(t_1, \dots,t_n)$ be a non-negative random variable depending on $t_1, \dots,
t_n$. Now choose $n$ points $t'_1, \dots, t'_n$ in $\R^d$ with respect to the truncated
distribution $\Psi'$ and define $Y'$ accordingly.  It is clear that

$$  \E(Y|A) =\E(Y'). $$

Let $c$ be a non-negative constant.  We say that $Y$ is {\it $c$-bounded } if
 $\E(Y|A) \le \Vol(B(R))^c$ and $\E(Y|B_i) \le \Vol (B(4^{i+1} R))^c$ for all $i \ge 0$.

 \begin{lemma} \label{lemma:Ybounded} If $Y$ is $c$-bounded then

 $$| \E(Y) - \E(Y') | = O(\E(Y) (\log n)^{-C_0 +cd/2}). $$
 \end{lemma}

 {\bf \noindent Proof of Lemma \ref{lemma:Ybounded}.} We start with the indentity

 $$\E(Y) = \E(Y|A)\P(A) + \E(Y |\overline{A})\P(\overline{A}). $$

  \noindent Since $\E(Y|A) =\E(Y')$, the triangle inequality implies that

  \begin{equation} \label{equ:Ystrong1} | \E(Y) - \E(Y') |
  \le \E(Y') \P(\overline{A})+\E(Y|\overline{A})\P(\overline{A}). \end{equation}

  \noindent To estimate $\E(Y|\overline{A})$, observe that

  \begin{equation} \label{equ:Ystrong2} \E(Y|\overline{A}) = \sum_{i=0}^{\infty} \E (Y|B_i \overline{A}) \P(B_i|\overline{A}). \end{equation}

  \noindent The ($c$-boundedness) assumption of the lemma implies

  $$ \E (Y|B_i \overline{A}) =  \E (Y|B_i) \le \Vol (B(4^{i+1} R))^c = O( 4^{cd(i+1)} R^{cd}) =
  O(4^{cd(i+1)} (\log n)^{cd/2}). $$

  \noindent Furthermore, as $B_i$ implies $\overline{A}$,

  $$\P(B_i|\overline{A}) = \frac{\P(B_i)}{ \P(\overline{A})} = O((\log n)^{C_0} \P(B_i)). $$

  \noindent On the other hand, $\P(B_i)$ is at most the probability that there is a point
  outside $B(4^iR)$. By the union bound and \eqref{equ:contentofball}, this probability is

  \begin{equation} \label{equ:Ystrong3}
  O( n \Psi (\overline {B(4^i R)}) = O( n \exp(- 4^{2i} R^2/2) (4^i R)^{d-2}).
  \end{equation}

\noindent For $i=0$, the right hand side of \eqref{equ:Ystrong3} is $\Theta ((\log
n)^{C_0})$ by the definition of $R$. For $i \ge 1$, the right hand side  of
\eqref{equ:Ystrong3}  is at most $n^{-2i}$, as

$$\exp(- 4^{2i} R^2/2)  = n^{(-1+o(1)) 4^{2i}} \le n^{-2i-1}. $$

 \noindent This shows that

$$\sum_{i=0}^{\infty} \E (Y|B_i \overline{A}) \P(B_i|\overline{A})
 = O\Big( \sum_{i=0}^{\infty} 4^{cd(i+1)} (\log n)^{cd/2} n^{-2i}  \Big)= O((\log n)^{cd/2}). $$

 \noindent Therefore the right hand side of \eqref{equ:Ystrong1} is at most

 $$O((\log n)^{cd/2} ) \P(\overline{A})= O((\log n)^{-C_0 + cd/2}), $$

 \noindent  proving the lemma. \hs

 Let $Y$ be the volume.
 It is clear that $Y$ is $1$-bounded.
 Applying Lemma \ref{lemma:Ybounded}, we have

 $$| \E\Vol (K'_n) - \E\Vol (K_n)| = O(\E \Vol (K_n)(\log n)^{-C_0 +d/2} )=O((\log n)^{-C_0 +d} ), $$

 \noindent since $\E\Vol (K_n) = \Theta ((\log n)^{d/2})$.  Moreover
 $\Var \Vol (K_n) = \Theta ((\log n)^{(d-3)/2})$. By setting $c_0$ sufficiently large,
 it  thus follows that

 $$| \E\Vol (K'_n) - \E\Vol (K_n)| = O( \sqrt {\Var \Vol (K_n) } (\log n)^{-C_0/2} ). $$

 We will use this estimate for proving the statement about the difference between the two
 variances. But first, let $Y$ be the square of the volume. It is clear that $Y$ is $2$-bounded. Thus, Lemma
 \ref{lemma:Ybounded} yields

 $$| \E\Vol (K'_n)^2 - \E\Vol (K_n)^2| = O(\E(\Vol (K_n))^2(\log n)^{-C_0 +d})= O((\log n)^{-C_0 +3d}), $$

 \noindent since $\Vol (K_n)^2 =O((\log n)^ {2d})$, which (by the definition of variance) implies,

 $$| \Var\Vol (K'_n) - \Var\Vol (K_n)| = O((\log n)^{-C_0 +3d } ) +
 |(\E\Vol (K'_n))^2- (\E\Vol (K_n))^2 | .$$

 \noindent On the other hand,

 $$ |(\E\Vol (K'_n))^2- (\E\Vol (K_n))^2 | = |\E\Vol (K'_n) +  \E\Vol (K_n) | |\E\Vol (K'_n)- \E\Vol (K_n) |, $$

 \noindent where $ |\E\Vol (K'_n)- \E\Vol (K_n) |$ is
 $O((\log n)^{-C_0 +d} )$ by the previous argument. Furthermore

 $$ |\E\Vol (K'_n) +  \E\Vol (K_n) | = O( \E\Vol (K_n)) = O((\log n)^{d/2}). $$

 \noindent Putting everything together, we obtain
 \begin{eqnarray*}
 | \Var\Vol (K'_n) - \Var\Vol (K_n)| &=& O ((\log n)^{-C_0+ 3d})
 + O( (\log n)^{-C_0 +d}(\log n)^{d/2})\\& =& O( (\log n)^{-C_0 + 3d}).
 \end{eqnarray*}

\noindent Again, by setting $c_0$ large, we have

 $$| \Var\Vol (K'_n) - \Var\Vol (K_n)| = O(  {\Var \Vol (K_n) } (\log n)^{-C_0/2} ), $$

 \noindent as claimed.

 To bound the difference between the two probabilities, define

 $$Y = \BI_{ \Vol (K_n) \ge t}. $$

 \noindent In this case, $Y$ is bounded from above by 1, thus it is $0$-bounded.
 Since $\E(Y) = \P( \Vol (K_n) \ge t)$, the claim follows instantly.
 \hs

We have the following
\begin{corollary}\label{cor:varK'} $\var \Vol (K'_n) = \Theta \left((\log n)^{(d-3)/2}\right).$
\end{corollary}

\section{The second coupling} \label{section:couple2}

In this section we will show that the first three conditions of
Lemma~\ref{lemma:twodistribution} are satisfied for the random variables $\Vol (\Pi_n)$
and $\Vol (K'_n)$. The fourth condition is just Theorem~\ref{theo:Poisson1}, whose proof
will come later. The first three conditions of Lemma~\ref{lemma:twodistribution} are
stated next.

\begin{lemma} \label{lemma:coupling2-1} For all sufficiently large $n$ we have

$$ |\E \Vol (\Pi_n) -\E \Vol (K'_n)| \le n^{-1/2+o(1)}  \sqrt{\Var \Vol (K'_n)} $$

$$  |\Var \Vol (\Pi_n) -\Var  \Vol (K'_n)| \le n^{-1/2+o(1)} {\Var \Vol (K'_n)},  $$

moreover, the following holds for all $t$

$$ | \P( \Vol (\Pi_n) \le t)- \P(\Vol (K'_n)\leq t)| \le n^{-1/2+o(1)}. $$

\end{lemma}

This lemma plus Theorem~\ref{theo:Poisson1} imply Theorem~\ref{theo:CLT1}, that is, the
central limit theorem for $\Vol (K'_n)$, which, in turn, implies Theorem~\ref{theo:CLT}.
So we will still have to prove Theorem~\ref{theo:Poisson1}, a major task which is the
content of the next four Sections. We mention further that Lemma \ref{lemma:coupling2-1}
implies the following.

\begin{corollary}\label{cor:varPi}. $\var \Vol (\Pi_n) = \Theta \left((\log n)^{(d-3)/2}\right).$
\end{corollary}

\begin{remark} \label{remark:rateofconv3} Let us notice that when applying
Lemma \ref{lemma:twodistribution}, the dominating error term comes from Theorem
\ref{theo:Poisson1}. Indeed, the error terms come from the first coupling are at most
$(\log n)^{-C}$, where $C$ can be arbitrarily large. The error terms from Lemma
\ref{lemma:coupling2-1} is even smaller, $n^{-1/2+o(1)}$. This implies the estimate on
the error term in Remark \ref{remark:rateofconv0}.
\end{remark}

Lemma \ref{lemma:coupling2-1} is a consequence of the following lemma.

\begin{lemma} \label{lemma:coupling2-2}  Let $A$ be a constant at least 10.
For any integer $n'$  between $n$ and $n +A \sqrt{n \log n}$

$$ |\E \Vol (K'_{n'}) -\E \Vol (K'_n)| \le  n^{-1/2+o(1)}  $$

$$  |\Var \Vol (K'_{n'}) -\Var  \Vol (K'_n)|  \le n^{-1/2+ o(1)} . $$

Moreover, for all $t$,

$$| \P( \Vol (K'_{n'})\le t)-\P( \Vol (K'_n)\le t)| \le n^{-1/2+ o(1)} . $$
\end{lemma}

{\bf \noindent Proof of Lemma  \ref{lemma:coupling2-1} via Lemma
\ref{lemma:coupling2-2}.}
 Let $A$ be a constant at least 10. We will use the fact that
 the probability that a Poisson variable with mean $n$ falls outside the interval
 $I= [n -A \sqrt{n \log n}, n + A \sqrt{n \log n}]$ is less than $n^{-A/4}$.
 As $\Vol (\Pi_n)$ is bounded from above by $\Vol (B(R))$,  we have

$$\E \Vol(\Pi_n)= \sum_{n' \in I} \E(\Vol (K'_{n'})) \P(n =n') + O(n^{-A/4} \Vol (B(R)) .$$

As $\Vol (B(R)) = O((\log n ^{d/2}))$, the last term on the right hand side is
$O(n^{-A/4+o(1)}) =O(n^{-1})$ as $A \ge 10$. So the first statement of Lemma
\ref{lemma:coupling2-2} implies
\begin{eqnarray*}
|\E \Vol (\Pi_n)-\E \Vol (K'_n)|& \le& \sum_{n' \in I} |\E(\Vol (K'_{n'})) -\E (\Vol
(K'_n))|\P(n=n') + O(n^{-1})\\& \le& n^{-1/2+o(1)}.
\end{eqnarray*}
Taking into account the fact that $\E(\Vol(K'_n)) = \Theta ((\log n)^{d/2})$ and $\Var
(\Vol (K_n)) =\Theta ((\log n)^{(d-3)/2})$, one can deduce the first statement of Lemma
\ref{lemma:coupling2-1}. The third statement of the same lemma can be proved the same
way.

Now we turn to the second statement. For every number ${n'} $ in the interval $I$, let
$E_{n'} $ denote the event that $n'$ is sampled (according to the Poisson distribution
with mean $n$) and $E_0$ denote the event that the sampled number does not belong to the
interval. The events $E_{n'}$ (with $n' \in I$ or $n'=0$) form a partition of the space.
Thus,

$$\Var \Vol(\Pi_n) = \E _{n'}  (\Var (\Vol(\Pi_n)|E_{n'} )) + \Var \E(\Vol
(\Pi_n|E_{n'} ), $$

\noindent where $n' \in I$ or $n'=0$. Notice that $\Vol (\Pi_n)|E_{n'} = \Vol (K'_{n'})$.
The rest of the proof is a calculation similar to the one above and is left as an
exercise. \hs

Let $H(r)$ be a halfspace at distance $r>0$ from the origin. Define $r$ so that the
probability content of $H(r) \cap B(R)$ is $\gamma \log n/n$ for some large constant
$\gamma$. As $\Psi'(H(r))=\Theta (e^{-r^2/2}r^{-1})$, $r=\Theta(\sqrt {\log n})$. For the
proof of Lemma \ref{lemma:coupling2-2} we need the the following claim.

\begin{claim}\label{cl:bbbr} The constant $\gamma$ can be chosen so  that $K'_n$ contains
$B(r)$ with probability at least $1-\frac 1n$.
\end{claim}

We explain the proof of this claim after the proof of Lemma \ref{lemma:sandwich} in the
next section.

{\noindent \bf Proof of Lemma  \ref{lemma:coupling2-2}.} Let us consider a number   $n'$
as in the lemma.  Let $\Omega$ denote the product space $B(R)^n$, equipped with the
$n$-fold product of $\Psi'$. A point $P$ in $\Omega$ is an ordered set $(x_1,...,x_n)$ of
$n$ random points (we generate the points one by one).  The $x_i$ are the coordinates of
$P$. We use $Y(P)$ to denote the volume of the convex hull of $P$ and $\mu$ to denote the
expectation of $Y(P)$.

\begin{remark} $Y(P)$ is, of course,  just another way to express
$\Vol(K'_n)$. It is however more convenient to use this notation in the proof below as it
emphasizes the fact that $Y$ is a function from $\Omega$ to $\BBR$.
\end{remark}

Define $\Omega', P', \mu'$ similarly (with respect to $n'$). Let us first consider the
expectations. Consider a point $P' =(x_1, \dots, x_{n'})$ in $\Omega'$ and  the canonical
decomposition

$$P' = P \cup Q $$

\noindent where $P=(x_1, \dots, x_n)$ and $Q=(x_{n+1}, \dots, x_{n'})$.  In order to
compare $\mu$ and $\mu'$, we rewrite $\mu$ as

$$\mu  = \int_{\Omega'}  Y(P)   dP' .$$

\noindent We have

$$\mu'- \mu = \int_{\Omega'} \left(Y(P')-Y(P)\right) d P'. $$

Now we are going to decompose $\Omega'$ into three parts $\Omega'_1, \Omega'_2,
\Omega'_3$ as follows

\begin{itemize}

\item $\Omega'_1= \{ P'|  \Conv (P) \,\,\, \hbox{does not contain the ball }  \,\,\, B(r) \} $.

\item $\Omega'_2 = \{P'| \Conv (P) \,\,\, \hbox{ contains the ball }  \,\,\, B(r) \,\,\,
\hbox{and} \,\,\,  B(r) \,\, \hbox{does not contain} \,\,\, Q \}. $

\item $\Omega'_3 = \Omega' \backslash (\Omega'_1 \cup \Omega'_2) $.
\end{itemize}

The measure of $\Omega'_1$ is the probability that the convex hull of a set of $n$ random
points does not contain $B(r)$, which is $O( 1/n)$, according to Claim \ref{cl:bbbr}. The
measure of $\Omega'_2$ is bounded from above by the probability that $B(r)$ does not
contain $Q$. This probability, by the union bound, is at most

$$ |Q| \times \Psi'( \overline{B(r)}) = O( \sqrt{n \log n})  \times \frac{(\log n)^{O(1)} }{n} =
n^{-1/2+o(1)} . $$

Since $Y(P')$ and $Y(P)$ are at most the volume of $B(R)$, which is $O((\log n)^{d/2})$,
$Y(P')-Y(P)$ is $O((\log n)^d)$. Thus

\begin{equation} \label{equ:coupling2-1} \int_{\Omega'_1 \cup \Omega'_2}  \left(Y(P')-Y(P)\right) d P'
= O((\log n)^{d} n^{-1/2+o(1)}) = n^{-1/2+o(1)}. \end{equation}

To estimate the integral over $\Omega'_3$, recall that in this region, $\Conv(P) =
\Conv(P')$ since

$$P' \backslash P= Q \subset B(r) \subset \Conv(P). $$

\noindent It follows that

\begin{equation} \label{equ:coupling2-2} \int_{\Omega'_3}
 \left(Y(P') - Y(P)\right)\, d P'= 0. \end{equation}

\eqref{equ:coupling2-1} and \eqref{equ:coupling2-2} together imply that
\[
\mu' -\mu  =  n^{-1/2+o(1)} ,
\]
\noindent proving the first part of the lemma.

The third part of the lemma follows now directly: the measure of $\Omega'_1 \cup
\Omega'_2$ is at most $n^{-1/2+o(1)}$, and on the the rest of $\Omega'$ the polytopes
$\Conv P= K'_n$ and $\Conv P'=K'_{n'}$ coincide.

The proof for the variance is similar. Notice that the  variance of $\Vol(K_n)$ is

$$s= \int_{\Omega'} |Y(P)- \mu|^2 \,\,d P' $$

\noindent and  the variance of $\Vol(K_{n'})$ is

$$s'= \int_{\Omega'} |Y(P')- \mu'|^2 \,\, d P'. $$

\noindent We have
\begin{equation} \label{equa:ss'1} |s' -s| = |\int_{\Omega'}
\left((Y(P') -\mu')^2 - (Y(P)-\mu)^2 \right)\, d P'|  \le \int_{\Omega'} |\CD (P')| \,\,
d P'
\end{equation}

\noindent   where

$$\CD(P') = (Y(P')  -\mu')^2 - (Y(P)-\mu)^2 . $$

\noindent It is obvious that

$$\CD(P') =  \big((Y(P') -\mu') + (Y(P)-\mu) \big) \big((Y(P') -\mu') - (Y(P)
-\mu ) \big) . $$

\noindent By the triangle inequality,

$$|\CD(P')| \le (|Y(P)+Y(P')+ \mu +\mu')(|Y(P')-Y(P)| +|\mu'-\mu|). $$

Since $Y(P')$ and $Y(P)$ are at most the volume of $B(R)$, which is $O((\log n)^{d/2})$,
$|\CD|$ is $O((\log n)^d)$. Thus, by arguing as before,

\begin{equation} \label{equ:coupling2-3} \int_{\Omega'_1 \cup \Omega'_2}  |\CD(P')|\,\, d P' = O((\log n)^{d} n^{-1/2+o(1)}) = n^{-1/2+o(1)}. \end{equation}

To estimate the integral over $\Omega'_3$, notice that in this region, $\Conv(P) =
\Conv(P')$. Therefore,

$$ \int_{\Omega'_3 } |\CD(P')| \,\,  d P' \le \int_{\Omega'_3}
(|Y(P)+ Y(P')+\mu + \mu')  |\mu'-\mu| \,\, d P'. $$

 \noindent But we just proved that
 $|\mu'-\mu| \le n^{-1/2+o(1)}$. Furthermore,
 all $Y(P'), Y(P), \mu', \mu$ are bounded from above by the volume of $B(R)$, which is
 $O((\log n)^{d/2})$. So

\begin{equation} \label{equ:coupling2-4} \int_{\Omega'_3}
 \big((Y(P') -\mu') + (Y(P)-\mu) \big) |\mu-\mu'| dP' \le n^{-1/2+o(1)} . \end{equation}

\eqref{equ:coupling2-3} and \eqref{equ:coupling2-4} together imply that

\begin{equation}   \label{equ:coupling2-5}
|s'-s | \le n^{-1/2+o(1)},
\end{equation}

\noindent concluding the proof. \hs

\section{Sandwiching $K'_n$}\label{section:sandwich}

By definition, $K'_n$ is contained in $B(R)$. In this section we will show that $K'_n$
contains the ball $B(r)$ with high probability where the radius $r$ is very close to $R$.
Recall that $R$ is defined in (\ref{equ:defR}) via
\[ R^2=2\log n+ \log(\log n)^{c_0}.
\]

The definition of $r$ comes a little later, we set first $\rho>0$ via
\begin{equation}
\rho^2=2\log n- \log\log n+ \log (c\log \log n)^{-2} \label{eq:defrho}
\end{equation}
where $c$ is a constant to be specified soon. Choose a system of points $y_1,\dots,y_m$
from the sphere $S(\rho)$ maximal with respect to the property that, for $i\ne j$,
\[
|y_i-y_j| \geq 2c_1.
\]
As $\rho=\sqrt {2\log n} (1+o(1))$ as $n$ goes to infinity, we have, just as in
Claim~\ref{cl:size}
\[
 m =\Theta\left((\log n)^{\frac {d-1}2}\right).
\]
Define the halfspace $H_i^+=\{x \in \R^d:\; y_i\cdot x \geq \rho^2\}$ and the cap $C_i$
as
\[
C_i=H_i^+ \cap B\Big(\sqrt {\rho^2+c_1^2}\Big).
\]
These caps are pairwise disjoint, and for $x \in C_i$
\[
\psi(x)=\Theta\left(\frac {c\sqrt {\log n}\log\log n}{n}\right).
\]
As $\Vol C_i=\Theta((\log n)^{-1/2}$, we have
\begin{equation}\label{eq:psiC}
\Psi(C_i)=\Theta\left(\frac {c\log\log n}{n}\right) \mbox{ and
}\Psi'(C_i)=\Theta\left(\frac {c\log\log n}{n}\right),
\end{equation}
since $C_i \subset B(R)$.

Set now $r=\rho- 5c_1^2/\rho$; it is clear then that this $r$ satisfies
\begin{equation}\label{eq:defr}
5c_1^2 <\rho^2- r^2 < 10c_1^2.
\end{equation}

\begin{lemma} \label{lemma:sandwich} For every $C>0$ the constants $c,c_1$ can be chosen so that the following holds.
 $K'_n$ contains $B(r)$ with probability at least
$1- (\log n)^{-C}$.
\end{lemma}

\begin{remark} This lemma is an analogue of a result from \cite{BD} for the uniform model
(see Section 2 for the definition). It is also a similar to Lemma 4.2 from \cite{vvu1},
which was proved using VC-dimension techniques. While in those results the probability
that $K_n$ does not contain $B(r)$ is at most $n^{-C}$, here we have the weaker bound
$(\log n)^{-C}$. The same bound was required in the uniform model when $K$ is a polytope,
see \cite{br}.
\end{remark}

{\bf \noindent Proof.} We claim first that every halfspace $H(r)$ at distance $r$ from
the origin contains a $C_i$ for some $i=1,\dots,m$. Assume $y$ is the nearest point of
$H(r)$ to the origin. Then $|y|=r$ and $y^*=\rho y/r$ lies on $S(\rho)$. As the system
$y_1,\dots,y_m$ is maximal, there is a $y_i$ with $|y^*-y_i| <2c_1$. Define $\alpha \in
(0,\pi/2)$ by $\sin \alpha =c_1/\rho$; it follows that the angle between $y$ and any
vector from $C_i$ is at most $3\alpha$. Consequently, $C_i$ is contained in the halfspace
with normal $y$ and at distance $\rho \cos 3 \alpha$ from the origin. A simple
computation shows now that for large enough $n$
\[
\rho \cos 3 \alpha > \rho-\frac {5c_1^2}{\rho}=r.
\]

\begin{claim}\label{cl:Br} There is a constant $b>0$ depending only on $d$ such that for all large enough $n$
\[
\P(B(r) \setminus K'_n\ne \emptyset ) =O\left(\frac {(\log n)^{\frac {d-1}2}}{(\log
n)^{bc}}\right).
\]
\end{claim}

{\bf \noindent Proof.} If $B(r)$ is not part of $K'_n$, then there is a halfspace $H(r)$
at distance $r$ from the origin which is disjoint from the random sample $X_n$. Then
there is a cap $C_i \subset H(r)$. Then $C_i \cap X_n=\emptyset$. Consequently
\begin{eqnarray*}
\P(C_i \cap X_n &=&\emptyset \mbox{ for some }i)\leq \sum_{i=1}^m\P( C_i \cap
X_n=\emptyset)\\
 & \leq & \sum_{i=1}^m\left(1-\Psi'(C_i)\right)^n \leq
           m\left(1-b\frac {c\log\log n}{n}\right)^n \\
 & \leq & m \exp\{-bc \log \log n\} = \frac {m}{(\log n)^{bc}} =O\left( \frac {(\log n)^{\frac {d-1}2}}{(\log
 n)^{bc}}\right).
\end{eqnarray*}
Here $b$ is the constant coming from (\ref{eq:psiC}). \hs

Choosing the constants $c$ and $c_1$ suitably completes the proof. \hs

\begin{remark} It is the choice of $r$ from (\ref{eq:defrho}) and (\ref{eq:defr})
that produces the bound $(\log n)^{-C}$. Also this choice of $r$ gives the estimates in
the next section. For the CLT for the volume, we could have taken $\rho^2=2\log n -\log
(c\log n)^3$ and $r=\rho-5c_1^2/\rho$ as well. This would have given
\begin{equation}\label{eq:randR}
\Psi'(C_i)=\Theta\left(\frac {c' \log n}{n}\right),
\end{equation}
and $1/n^{-c'}$ for the probability that $K'_n$ does not contain $B(r)$. But this choice
does not work for $f_s(K_n)$ (see Remark 13.7). That's why we used (\ref{eq:defrho}) and
(\ref{eq:defr}) for the definition of $r$.
\end{remark}

The proof of Claim~\ref{cl:bbbr} goes along very similar lines. One can take $\rho
^2=2\log n -\log (\gamma' \log n)^3$, for instance, and use the same argument. We omit
the details.

One can prove similarly that $\Pi_n$ contains $B(r)$ with high probability. Here is the
quantitative statement, the routine proof is left to the interested reader.

\begin{lemma} \label{lemma:sandwichPi} For every $C>0$ the constants $c,c_1$ can
be chosen so that the following holds. $\Pi_n$ contains $B(r)$ with probability at least
$1- (\log n)^{-C}$.
\end{lemma}

\begin{remark} Note that $K_n$ is sandwiched between $B(R)$ and $R(r)$ with high probability, and both
$r,R=\sqrt {2 \log n}(1+o(1))$. This almost implies (\ref{eq:vol}) for the  expectation
of $\Vol (K_n)$, the only trouble being that $K_n$ can have arbitrarily large volume when
it is not contained in $B(R)$.
\end{remark}

\section{The dependency graph}  \label{section:depend}

With the notation of the previous section we define the annulus $A(R,r)=B(R) \setminus
B(r)$, and let $V_i$ denote the Voronoi region of $y_i$ ($i=1,\dots,m$). This means that
$x \in V_i$ if and only if $|x-y_i| \leq |x-y_j|$ for all $j$. The sets $W_i=V_i \cap
A(R,r)$ will be called {\sl cells} and will play an important role in the central limit
theorems. The following estimate will be needed.

\begin{claim}\label{cl:psiW} For each $i$
\[
\Psi'(W_i) = \Theta\left(\frac {\log \log n}{n}\right).
\]
\end{claim}

{\bf \noindent Proof.} This is quite simple and similar to (\ref{eq:psiC}) and is
therefore omitted. $\Box$

The {\sl dependency graph} $G(V,E)$ has, by definition, vertex set $V(G)=\{1,\dots,m\}$
and edge set $E(G)$ with $(i,j) \in E(G)$ if and only if there are $a_i \in W_i$ and $a_j
\in W_j$ and $b \in A(R,r)$ such that the segments $[a_i,b]$ and $[b,a_j]$ lie completely
in $A(R,r)$. In other words, if and only if $[a_i,b]\cap B(r)=\emptyset$ and $[a_j,b]\cap
B(r)=\emptyset$ for some $a_i \in W_i$, $a_j \in W_j$ and $b \in A(R,r)$. Let $D$ denote
the maximal degree in the dependency graph.

\begin{theorem}\label{th:maxdeg} $D =O\left((\log \log n)^{\frac {d-1}2}\right)$.
\end{theorem}

{\bf \noindent Proof.} This is a simple matter using elementary geometry. Observe first
that if the segment $[a,b] \subset A(R,r)$ and $2\gamma$ is the angle between vectors $a$
and $b$, then $\cos \gamma \ge r/R$. We can estimate $\sin \gamma$ using the definitions
of $R$ and $r$:
\begin{eqnarray*}
\sin \gamma &\leq & \sqrt{1 - \left(\frac rR\right)^2}= \frac 1R \sqrt {R^2-r^2}=\frac
1R \sqrt{R^2-\left(\rho -\frac {5c_1^2}{\rho}\right)^2} \\
   & \le & \frac 2R \sqrt{\log (\log n)^{2c_0}} =O\left(\sqrt{\frac {\log \log n}{\log n}}\right).
\end{eqnarray*}

Suppose next that $a_i \in W_i$ and let $2\alpha _i$ be the angle between $a_i$ and
$y_i$. Set $a_i^*=\rho a_i/|a_i| \in S(\rho)$. The maximality of the system
$y_1,\dots,y_m$ implies that $|a_i^*-y_i| \le 2c_1$, which, in turn, shows that $\sin
\alpha _i \le c_1/\rho$. Consequently $\alpha =O((\log n)^{-1/2})$.

Assume $(i,j) \in E(G)$ and let $a_i \in W_i$, $a_j \in W_j$ and $b \in A(R,r)$ be the
vectors such that the segments $[a_i,b]$ and $[a_j,b]$ are disjoint from $B(r)$. Let
$2\beta$ be the angle between vectors $y_i,y_j$. Then
\[
\beta \le \alpha_i+\gamma +\alpha _j=O\left(\sqrt{\frac {\log \log n}{\log n}}\right).
\]
This, of course, implies  that for $(i,j) \in E(G)$
\[
 |y_j-y_i| \leq 2R \sin \beta =O(\sqrt{\log \log n}).
\]
This means that all $y_j$ with $(i,j) \in E(G)$ are contained in a ball, centered at
$y_i$ and of radius $O(\sqrt{\log \log n})$. Since all $y_j \in S(\rho)$ and since they
are at distance $2c_1$ apart, the usual volume estimate gives the statement of the
theorem. \hs

We establish one more inequality here.

\begin{claim}\label{cl:volW} For each $i$
\[
\Vol (W_i) = \Theta \left(\frac {\log \log n}{\sqrt {\log n}} \right).
\]
\end{claim}

{\bf \noindent Proof.} For each $t \in [r,R]$, $W_i\cap S(t)$ has constant, that is,
$\Theta(1)$ $(d-1)$-dimensional volume, so $\Vol (W_i)=O(R-r)$, and
\[
R-r=\frac 1{R+r} (R^2-r^2)=\frac 1{R}\Theta (\log \log n)
\]
as we have seen in the previous proof. \hs

\section{Central limit theorem for the Poisson model} \label{section:CLTPoisson}

We are going to apply the Baldi-Rinott theorem for $\Pi_n$ conditioned on $B(r) \subset
\Pi_n$. This condition will be denoted by $B$. Recall from Lemma~\ref{lemma:sandwichPi}
that
\[
\P(B(r)\subset \Pi_n) \ge 1- (\log n)^{-C}.
\]

Assume condition $B$ holds and define the random variable $\xi_i=\Vol (W_i \cap \Pi_n)$.
Clearly, $\xi:=\sum_1^m \xi_i =\Vol (\Pi_n)- \Vol (B(r))$. This shows that, under
condition $B$, the CLT for $\xi$ holds if and only if it holds for $\Vol (\Pi_n)$.

\begin{claim}\label{cl:baldi} Assume condition $B$ holds. Given disjoint subsets $V_1,V_2$
of the vertex set of the dependency graph with no edge between them, the random variables
$\{\xi_i:\; i \in V_1\}$ are independent of the random variables $\{\xi_j:\; j \in
V_2\}$.
\end{claim}

{\bf \noindent Proof.} The intersection $W_i \cap \Pi_n$ is determined by the facets of
$\Pi_n$ intersecting $W_i$. These facets are determined by their vertices. If there are
no common vertices for the facets intersecting the $W_i$ with $i \in V_1$ and the $W_j$
with $j \in V_2$, then the corresponding $\xi_i$ are independent. This is exactly how the
dependency graph has been defined. \hs

Write $\P^*$, $\e^*$, $\var^*$ for $\P$, $\e$, $\var$ under condition $B$. In the next
section we will prove the following estimates.

\begin{lemma} \label{lemma:baldii} We have

$$ |\E^* \Vol (\Pi_n) -\E \Vol (\Pi_n)| \le (\log n)^{-C_0/4}  \sqrt{\Var \Vol (\Pi_n)}, $$

$$  |\Var^* \Vol (\Pi_n) -\Var  \Vol (\Pi_n)| \le (\log n)^{-C_0/4} {\Var \Vol (\Pi_n)}, $$

$$ | \P^*( \Vol (\Pi_n) \le t)- \P(\Vol (\Pi_n)\leq t)| \le (\log n)^{C_0/4}. $$
\end{lemma}

The inequality for the variances shows that

$$\var^* \Vol (\Pi_n) =\Omega (\var \Vol (\Pi_n)) =\Omega \left((\log n)^{\frac
{d-3}2}\right)$$.

We have seen that the maximal degree in $G$ is $O((\log \log n)^{(d-1)/2})$
(Theorem~\ref{th:maxdeg}), and $\xi_i= \Vol (W_i) = O(\log \log n/\sqrt {\log n})$. So
the Baldi-Rinott theorem applies and gives the following CLT.

\begin{theorem} \label{theo:Poisoncond}
Let $d$ be a fixed integer at least $2$.  For any value of $t$,

\[ | \P^* \Big(\frac{|\Vol (\Pi_n) -
\E^*\Vol (\Pi_n)|}{ \sqrt {\Var^* \Vol (\Pi_n)}} \le t \Big) -\Phi (t) | \\
 = O\Big(\frac {(\log \log n)^{\frac {d+4}2}}{(\log n)^{\frac {d-1}4}}\Big).
\]

\end{theorem}\hs

This theorem and Lemma \ref{lemma:baldii} show that $\Vol (\Pi_n)$ and $\Vol (\Pi_n)|B$
satisfy conditions of Lemma \ref{lemma:twodistribution}. So our main central limit
theorem, Theorem~\ref{theo:CLT}, follows as soon as we prove Lemma~\ref{lemma:baldii}.
This is our next (and final) task.

\section{Proof of Lemma \ref{lemma:baldii} } \label{section:final}

This is similar to, and much simpler than, the proof in Section~\ref{section:couple1}.
The first step is a copycat of Lemma \ref{lemma:coupling1}.

\begin{lemma} \label{lemma:coupling3} Let $B$ denote the condition that $B(r) \subset K'_n$.
Then we have, for large enough $n$,

$$| \E(\Vol (K'_n)|B)- \E\Vol (K'_n)| \le (\log n)^{-C_0/2} $$

$$| \Var (\Vol(K'_n)|B)- \Var \Vol(K'_n)| \le   (\log n)^{-C_0/2}. $$

\noindent Furthermore, for all $t$,

$$|\P (\Vol (K'_n) \le t|B) - \P (\Vol (K'_n) \le t) | \le (\log n)^{-C_0/2}. $$
\end{lemma}

{\bf \noindent Proof.} We use the first few lines of the proof of Lemma
\ref{lemma:Ybounded} with condition $A$ replaced by $B$, events $B_i$ do not appear yet.
Then (\ref{equ:Ystrong1}) says that
\begin{equation}\label{eq:bald1}
|\E(Y|B)-\E(Y)| \leq (\E(Y|B)+\E(Y|\overline{B}))\P(\overline{B}),
\end{equation}
where $Y=Y(t_1,\dots,t_n)$ is a $c$-bounded, nonnegative random variable.

When $Y$ is just the volume, $Y$ is bounded by $O((\log n)^{d/2})$ so its expectation,
under any condition, is bounded the same way. Since $\P(\overline{B}) \le (\log
n)^{-C_0}$ by Lemma \ref{lemma:sandwich}, we are finished with the first inequality.

The third is proved by setting $Y=\BI_{\Vol (K'_n)\le t}$. The second inequality follows
the same way as the corresponding inequality for variances in Lemma
\ref{lemma:coupling1}. \hs

 We show finally how this lemma implies Lemma \ref{lemma:baldii}.

{\bf \noindent Proof} of Lemma \ref{lemma:baldii}. We give the proof for $\E$ first. As
before, write $E_n'$ for the event that $|X(n)|=n'$.
\begin{eqnarray*}
|\E^* \Vol (\Pi_n) &-& \E \Vol (\Pi_n)|= \left|\sum _0^{\infty} \Big(\E(\Vol
(K'_{n'})|B)-\E
\Vol K'_{n'}\Big)\P(n=n')\right|\\
&\le & \sum_{n' \in I} (\log n')^{-C_0/2}\P(n=n') + O((\log n)^{d/2}n^{A/4})\\
& = & O((\log n)^{-C_0/2}).
\end{eqnarray*}
This suffices for the the expectations as $\Var \Vol (\Pi_n) =\Theta ((\log n)^{(d-3)/2}$
by Corollary \ref{cor:varPi}. Of course, we chose $C_0$ large enough.

The proof for $\var^*$ and $\P^*$ is similar and is left to the reader. \hs

We want to emphasize here that the proofs of Theorems \ref{theo:Poisson1},
\ref{theo:CLT1}, and \ref{theo:CLT} have finally been completed at this point.

\section{Proof of Theorem \ref{theo:CLTfaces} } \label{section:faces}

The proof of Theorem \ref{theo:CLTfaces} follows the plan in Section \ref{section:plan}
closely. In fact, most of the arguments are the same as in the proof of Theorem
\ref{theo:CLT}, except for a few technical modifications, and a single extra difficulty:
finding the right bound $M$ on the number of $s$-faces intersecting cell $W_i$. Thus,
instead of working out all details, we only state the main steps and point out what
modifications are needed, plus explain how the bound $M$ can be found.

We have seen in Theorem~\ref{th:varfs} that the variance satisfies

$$\Var (f_s(K_n))= \Theta ((\log n)^{(d-1)/2}). $$

\subsection {The first coupling} Lemma \ref{lemma:coupling1} still holds if one replaces
$\Vol$ by $f_s$. Notice that the proof of this lemma only requires the $c$-bounded
property. The number of faces has this property (for some sufficiently large constant
$c$).  Indeed, one can show that with very high probability (say $1- n^{-100d}$) the
number of vertices is at most $(\log n)^{d}$. This, together with a simple geometric
argument shows that the number of faces is $c$-bounded for some constant $c$. The same
proof goes for the square of the number of faces.

After the first coupling, it is left to prove the following variant of Theorem
\ref{theo:CLT1}.

\begin{theorem} \label{theo:CLT1faces} Let $s$ be an integer between $0$ and $d-1$.
There is a function $\ep(n)$ tending to zero as $n$ tends to infinity such that for all
$t$

$$ | \P \Big( \frac{f_s(K'_n)- \E f_s(K'_n)} {\sqrt{\Var f_s(K'_n)}}  \le t \Big)-\Phi (t) | \le \ep(n). $$

\end{theorem}

\subsection {The second coupling} The proof for the second coupling is almost the same as before.
A small technical modification one needs to make here is to introduce a new part
$\Omega'_0$ in the partition which contains those $P'$ where $\Conv (P')$ has more than
(say) $(\log n)^d$ vertices. The probability of $\Omega'_0$ will be less than $n^{-1/2}$.
Now define $\Omega'_3 =\Omega \backslash (\Omega'_0 \cup \Omega'_1 \cup \Omega'_2)$. The
rest of the proof is the same. In fact, since both the expectation and variance of
$f_s(K'_n)$ are also polylogrithmic in $n$ (similar to those of the volume), the error
term $n^{-1/2+o(1)}$ remains unchanged in all these estimates.

After the second coupling one needs the $f_s$ variant of Theorem \ref{theo:Poisson1}.

\begin{theorem} \label{theo:Poisson1faces}
Let $d$ be a fixed integer at least $2$ and $0 \le s \le d-1$.  There
is a function $\ep(n)$ tending to $0$ as $n$ tends to infinity such that the following holds.
For any value of $t$,

\begin{equation} \label{equ:CLTPoissonf} | \P \Big(\frac{|f_s (\Pi_n) -
\E f_s (\Pi_n)| }{ \sqrt {\Var f_s (\Pi_n)}} \le t \Big) -\Phi (t) |
\le \ep(n).
\end{equation}

\end{theorem}

\begin{remark} \label{remark:rateofconv2faces} One can take $\ep(n)= (\log n)^{-(d-1)/4+o(1)}$. This error
term will be the dominating one when we apply, twice, Lemma \ref{lemma:twodistribution} .
\end{remark}

\subsection {The dependency graph} The dependency graph is the same as before with
\[
m=\Theta ((\log n)^{(d-1)/2}), \;\;\; D=O((\log \log n)^{(d-1)/2}),\mbox{ and }\Psi'
(W_i)=\Theta((\log \log n)/n).
\]

For proper accounting $f_s(\Pi_n)$ we have to define the random variable $\xi_i=f(W_i,s)$
suitably. Fotr this purpose we use Reitzner's method from \cite{rei1}. For an
$s$-dimensional face, $L$, of $\Pi_n$, let $f(W_i,L)$ denote the number of vertices of
$L$ contained in $W_i$, and set
\[
f(W_i,s)= \frac 1{s+1}\sum_L f(W_i,L).
\]
Since $\Pi_n$ is simplicial and has no vertex on the boundary of any $W_i$ with
probability one, $f_s(\Pi_n)=\sum _{i=1}^m f(W_i,s)$. The expected number of $|X(n) \cap
W_i|= \Theta(\log \log n)$, which, in turn, shows that that the expectation of $f(W_i,s)$
is $\Omega(\log \log n)$. But there is an extra difficulty here: we need a bound $M$ on
each $f(W_i,s)$ when applying the Baldi-Rinott theorem. The condition $B(r) \subset
\Pi_n$ is not enough and we have to introduce a new condition, to be denoted by $B_i$:
\[
|X(n) \cap W_i| \leq c_2\log \log n \mbox{ for each } i.
\]
where $c_2$ is a large constant. It is straightforward to check that for any $C>0$, $c_2$
can be chosen so large that
\[
\P(B_i \mbox{ holds}) \ge 1-(\log n)^{-C}.
\]
Then the union bound shows that
\[
\P(B_i \mbox{ fails for some } i) =O((\log n)^{-C+(d-1)/2}).
\]

It is clear that if $L$ is an $s$-face of $\Pi_n$ contributing to $F(W_i,s)$, then all
vertices of $L$ belong to a cell $W_j$ with $i,j$ connected in $G$ or to $W_i$. There are
at most $D$ such cells. So under condition $B_i$, there are at most $c_2D\log \log n$
vertices in the union of these cells. This shows that $M= (\log \log n)^{d^2}$ works and
the application of the Baldi-Rinott theorem goes through.

Again we have to remove the conditions $B,B_1,\dots,B_m$. This is done in the same way as
in Section~\ref{section:final}.

\begin{remark} This is where the careful choice of $r$ (in fact, $\rho$) pays off. With
the more generous selection $\rho^2=2\log n-\log(c \log n)^3$, we would only have
$f(W_i,s)=O((\log n)^{d/2})$, and the right hand side in the estimate of the Baldi-Rinott
theorem does not tend to zero.
\end{remark}

\section{Concluding remarks} \label{section:remarks}

Our plan can be used for many other parameters. In certain cases, one merely has to repeat the proof.
In others, however, there are substantial technical difficulties. Let us present two representative examples.

\vskip2mm

{\it The surface area of $K_n$.} The proof is more or less  the same as the  proof for
the volume. The reader is invited to work out the details. In fact, the result holds for
all intrinsic volumes, but the estimate for variance is not straightforward.

\vskip2mm

{\it The probability content of $K_n$.} The probability content of $K_n$ is $\Psi(K_n)$.
For this parameter, the general plan still works, but there is a non-negligible
difficulty. In the proof of the second coupling, we used the fact that the expectation
and variance of the random variable under study (such as the volume, number of faces, or
even the surface area)  are both polylogarithmic in $n$. Thus, the error term
$n^{-1/2+o(1)}$ is dominating and one can finish the proof easily. For the case of the
probability content, it is no longer true, as the variance is  $n^{-2 +o(1)}$. To
overcome this obstacle, we can follow \cite{vvu2} and start by proving a sharp
concentration result, which gives a tight control on the tail $Y(P)-\mu$ and
$Y(P')-\mu'$.   Such a concentration result is available thanks to the method developed
in \cite{vvu1}.  The details will appear elsewhere.

\end{document}